\documentclass[10pt,a4paper]{article}
\makeatletter
\usepackage[left=3cm,top=3cm,bottom=3cm,right=3cm,head=1cm,foot=1cm]{geometry}
\usepackage{amssymb}
\usepackage{amsmath}
\usepackage{amsthm}
\usepackage{amsfonts}
\usepackage{amssymb}
\usepackage{amscd}
\usepackage{url}
\usepackage[all]{xy}

\newcommand{\h}{\mathbf{H}}

\newcommand{\g}{\mathfrak{g}}

\newcommand{\5}{\hspace{0,5cm}}

\newcommand{\ad}{\stackrel{\mathrm{ad}}{\triangleright}}

\newcommand{\bad}{\stackrel{\mathrm{ad}}{\blacktriangleright}}
\title{On Dirac Operators and Spectral Geometry of Compact Quantum Groups}
\author{Antti J. Harju}
\date{}
\begin{document}
\maketitle

\section*{Introduction}

The classical Dirac operator $D$ on a Lie group $G$ with a Lie algebra $\g$ can be seen as a purely algebraic object living in the noncommutative Weyl algebra $U(\g) \otimes \mathrm{cl}(\g)$, see \cite{HP}. It is a self adjoint element and $D^2$ is a sum of Casimir elements in $U(\g) \otimes \mathrm{cl}(\g)$. Furthermore, $D$ is equivariant in the sense that there exists a Lie algebra homomorphism $\g \rightarrow U(\g) \otimes \mathrm{cl}(\g)$ and $D$ commutes with its image.

The algebra $U(\g) \otimes \mathrm{cl}(\g)$ acts on the Hilbert space $L^2(G) \otimes \Sigma$, where $\Sigma$ is an irreducible $\mathrm{cl}(\g)$-module making $D$ an unbounded Fredholm operator with infinitely many positive and negative eigenvalues. The Dirac operator has an important role in index theory and $K$-homology. In the case of compact group $G$ the Dirac operator is a fundamental object in the spectral geometry: the spectral triple $(C^{\infty}(G), D, L^2(G) \otimes \Sigma)$ defines an alternative operator theoretic approach to the Riemannian geometry. 

The first attempt to define a Dirac type operator in deformed settings, more precisely for $SU_q(2)$, was made in \cite{BK}. This approach provides a deformation of the Dirac operator so that its algebraic properties survive. Especially it defines an equivariant system. As a Hilbert space operator it is unbounded with infinite number of positive and negative eigenvalues and therefore provides an interesting application in the index theory. However, the spectrum of this operator grows exponentially and therefore it cannot be applied in a spectral triple. A generalization of this approach was done in \cite{H} where the construction was given for any quantum group based on a simple Lie algebra.

An isospectral deformation of the spectral triple on $SU_q(2)$ was done in \cite{DV}, \cite{CP}. In this approach one defines a natural Hilbert space for the Dirac operator and decomposes it into irreducible components for the action of the symmetry algebra $U_q(\mathfrak{su_2})$. Then the Dirac operator is given a constant action on each component. In order to make the operator satisfy the properties of a spectral triple its spectrum is chosen to match with the spectrum of a classical Dirac operator. Such a system is automatically equivariant. In the approach \cite{CP} this type of an operator was defined on $L^2(G_q)$ which is a Hilbert space completion of the algebra of polynomial functions on $SU_q(2)$. This is not really a deformation because the classical Dirac operator does not act on the space of $L^2$ functions, however, the operator fits into a definition of a spectral triple and therefore provides a well defined noncommutative space. In \cite{DV} similar operator was constructed on a Hilbert space $L^2(G_q) \otimes \mathbb{C}^2$ which notices the spinor module $\mathbb{C}^2$. Again this leads to a well defined spectral triple. The relationship between these two approaches was found in \cite{CP2}. 

In \cite{NT} an equivariant Dirac opeator was defined for any quantum group deformation of a simple, compact and simply connected Lie group which satisfies the axioms of a spectral triple. The model is based on conjugating the classical Dirac operator with a unitary twist $F$. It is known that the construction is independent of any choices, such as a twist, up to unitary equivalence of the spectral triples, \cite{NT2}. As $C^*$-algebras, the quantum groups are $KK$-equivalent with their classical limit \cite{N}. The $K$-homology cycles defined by \cite{NT} corresponds to the fundamental $K$-homology cycle defined by the classical Dirac operator under this equivalence \cite{NT2}. 

Here we give a short review of the Dirac operators in the approaches \cite{H} and \cite{NT}. We construct several examples of these Dirac operators and use those to build Fredholm modules and spectral triples. We consider the case $SU_q(2)$ with details and find how the Dirac operator in the algebraic approach \cite{H, BK} and in the geometric approach \cite{NT} are related. Furthermore, we see that the isospectral deformation of \cite{DV} is essentially the same as the spectral triple of \cite{NT} applied in the case $SU_q(2)$. We shall use the details of $SU_q(2)$ to create a spectral triple for the quantum group $U_q(2)$ and study its properties. We study the quantum sphere $S^2_q$ which is by definition a fixed point algebra $SU_q(2)^{U(1)}$. It follows that the algebraic approach \cite{H} leads to a well defined $0$-summable spectral triple. This model was found earlier in a case study \cite{S}. We also describe how to construct the operator in the algebraic approach and write down formulas in the case $SU_q(3)$. 

\section{Quantum Group Preliminaries}

\textbf{1.1.} Denote by $G$ a simple, simply connected and compact Lie group and $\g$ its Lie algebra. Consider a fixed maximal torus and let $\mathfrak{h} \subset \g$ be a Cartan subalgebra and $\{ \alpha_i: 1 \leq i \leq n \}$ a set of simple roots. Let $(a_{ij}): 1 \leq i \leq n$ denote a Cartan matrix of $\g$ and $\{d_i: 1 \leq i \leq n\}$ coprime positive integers such that $(d_{i} a_{ij})$ is a symmetric matrix. $P_+$ is the set of dominant integral weights of $G$. $(\cdot, \cdot)$ denotes the Killing form if not otherwise specified.

Fix a finite dimensional irreducible unitary representation ($V_{\lambda}, \pi_{\lambda}$) of $\g$ for each $\lambda \in P_+$. Denote by $W^*(G)$ the Hopf von Neumann $*$-algebra generated by the fixed representations of $\g$, i.e., $W^*(G)$ is the $l^{\infty}$-direct sum of $B(V_{\lambda})$. Define  $*$-algebra $U(G) = \prod_{\lambda} B(V_{\lambda})$ of unbounded densely defined operators affiliated with $W^*(G)$. The primitive coproduct of $W^*(G)$ extends to a $*$-homomorphism $\triangle: U(G) \rightarrow \prod_{\lambda,\nu} B(V_{\lambda} \otimes V_{\nu}) :=U(G \times G)$.

Let $q \in (0,1)$. Define an associative noncocommutative Hopf $*$-algebra $U_q(\g)$ which is the polynomial algebra generated by $e_i, f_i, k_i, k_i^{-1}: 1 \leq i \leq n$, subject to the relations 
\begin{eqnarray*}
& &[k_{i}, k_{j}] = 0,\5  k_i k_i^{-1} = 1 \5 k_{i} e_j k_{i}^{-1} = q_i^{a_{ij}/2} e_j,\5 k_{i} f_j k_{i}^{-1} = q_i^{-a_{ij}/2} f_j, \\
& &[e_i, f_j] = \delta_{ij}\frac{k^2_{i} - k_{i}^{-2}}{q_i -q_i^{-1}}, \5 q_i = q^{d_i}
\end{eqnarray*}
and the $q$-Serre relations \cite{Dr,Ji}. We choose the Hopf $*$-structure 
\begin{eqnarray*}
& &\triangle_q(k_i) =k_i \otimes k_i, \5 \triangle_q(e_i) = e_i \otimes k_i + k_i^{-1} \otimes e_i,\5 \triangle_q(f_i) = f_i \otimes k_i + k_i^{-1} \otimes f_i, \\
& &S_q(e_i) = -q e_i,\5 S_q(f_i) = -q^{-1}f_i,\5 S_q(k_i) = k_i^{-1}, \\
& & \epsilon_q(k_i) = 1,\5 \epsilon_q(e_i) =  \epsilon_q(f_i) = 0, \5 e_i^* = f_i, \5 f_i^* = e_i,\5 k_i^* = k_i.
\end{eqnarray*}
The equivalence classes of irreducible finite dimensional representation of $U_q(\g)$ are classified by the integral dominant weights of $\g$. A module of highest weight $\lambda $ with a highest weight vector $|\lambda, \lambda \rangle$ has the defining properties 
\begin{eqnarray*}
\pi_{q,\lambda}(k_i) |\lambda, \lambda \rangle = q_i^{\lambda(h_i)/2} |\lambda, \lambda \rangle, \5 \pi_{q,\lambda}(f_i) |\lambda, \lambda \rangle = 0
\end{eqnarray*} 
for each $i$, where $h_i \in \mathfrak{h}$ so that $\alpha_j(h_i) = a_{ij}$. The dimensions of each weight spaces of a representation of $U_q(\g)$ are the same as for the representation of $\g$, \cite{Lu, Ro}.

Let us fix an irreducible unitary representation $(V_{q,\lambda}, \pi_{q,\lambda})$ of $U_q(\g)$ for each $\lambda \in P_+$. Define $W^*(G_q)$ the Hopf von Neumann $*$-algebra generated by the fixed representations of $U_q(\g)$ and define $*$-algebra $U(G_q) = \prod_{\lambda} B(V_{\lambda,q})$ of unbounded densely defined operators affiliated with $W^*(G_q)$. The comultiplication of $W(G_q)$ extends to a $*$-homomorphism on $\triangle_q: U(G_q) \rightarrow \prod_{\lambda, \nu} B(V_{q,\lambda} \otimes V_{q,\nu}):= U(G_q \times G_q)$. The $*$-algebra $U(G_q)$ is thought as a completion of $U_q(\g)$. 

The noncocommutativity of the coproduct is controlled by the $R$-matrix, $R \in U(G_q \times G_q)$ so that 
\begin{eqnarray*}
\sigma(\triangle_q(x)) = R \triangle_q(x) R^{-1} 
\end{eqnarray*}
where $\sigma$ is the flip automorphism.

For each $\lambda \in P_+$ there is a $*$-algebra isomorphism $\phi_{\lambda}: B(V_{\lambda,q}) \rightarrow B(V_{\lambda})$ which identifies the centers. We can apply these isomorphism to define a $*$-algebra isomorphism
\begin{eqnarray*}
\phi: W^*(G_q) \rightarrow W^*(G) 
\end{eqnarray*}
which identifies the centers and extends to a $*$-isomorphism $\phi: U(G_q) \rightarrow U(G)$. However, the coproducts of $U(G)$ and $U(G_q)$ do not respect the isomorphism. It is proved in \cite{NT} that there exists a $*$-isomorphism $\phi: W^*(G_q) \rightarrow W^*(G)$ identifying the centers and a unitary $F\in W^*(G) \otimes W^*(G)$ such that
\begin{eqnarray*}
&& (\phi \otimes \phi) \triangle_q(x) = F \triangle(\phi(x)) F^{*}, \5 \hbox{for all $x \in W^*(G_q)$}, \\
&& (\epsilon \otimes \mathrm{id}) F = ( \mathrm{id} \otimes \epsilon) F  = 1, \nonumber \\
&& (\phi \otimes \phi) R = F_{21} q^{- \sum_k x_k \otimes x_k}  F^{*}, \nonumber \\
&& \hbox{the associator $\Phi = (\mathrm{id} \otimes \triangle)(F^*)(1 \otimes F^*)(F \otimes 1)(\triangle \otimes 1)(F)$ coincides with $\Phi_{\mathrm{KZ}}$}, \nonumber
\end{eqnarray*}
where $\{x_k \}$ is the basis of $\g$ fixed from the condition $(x_k, x_l) = - \delta_{kl}$ and the Drinfeld associator $\Phi_{\mathrm{KZ}}$ is the associativity morphism of tensor products in the category of Lie algebra representations which is determined by the monodromy of Knizhnik-Zamolodchikov equations. We shall assume that $F$ has these properties in the following. The twist is not unique. \\

\textbf{1.2.} Consider the dual Hopf algebra $\mathbb{C}[G_q]$ of $U_q(\g)$ which reduces to the Hopf algebra of regular (or representative) functions on $G$ for $q \rightarrow 1$. As a vector space $\mathbb{C}[G_q]$ is spanned by the matrix elements of the endomorphisms of irreducible finite dimensional module $V_{\lambda}$
\begin{eqnarray*}
\mathbb{C}[G_q] = \bigoplus_{\lambda \in P_+} V^*_{q,\lambda} \otimes V_{q,\lambda}.
\end{eqnarray*}
Let us equip each module ($V_{q,\lambda}$, $\pi_{q,\lambda})$ with an inner product and fix an orthonormal basis $\{ |\lambda, \nu \rangle : \nu \in I\}$ ($I$ is an index set) for each $V_{q,\lambda}$. The algebra $U_q(\g)$ acts on a basis vector  $t^{\lambda}_{\mu, \nu} = \langle \lambda, \mu | \otimes |\lambda, \nu \rangle$ of $\mathbb{C}[G_q]$ from the left by 
\begin{eqnarray*}
\partial(x) t^{\lambda}_{ \mu, \nu }  =  \langle \lambda, \mu | \otimes ( \pi_{q,\lambda}(x)| \lambda, \nu \rangle )
\end{eqnarray*}
for all $x \in U_q(\g)$. 
This extends to an action of $W^*(G_q)$ and $U(G_q)$ on $\mathbb{C}[G_q]$. The pairing $ \mathbb{C}[G_q] \otimes U_q(\g) \rightarrow \mathbb{C}$ is defined by
\begin{eqnarray*}
 t^{\lambda}_{\mu \nu}(x) = \langle \lambda, \mu | ( \pi_{q,\lambda}(x)| \lambda, \nu \rangle )
\end{eqnarray*} 
for all $x \in U_q(\g)$. 

The space $\mathbb{C}[G_q]$ can be equipped with a multiplication fixed from the condition 
\begin{eqnarray*}
(t^{\lambda}_{\mu,\nu} t^{\lambda'}_{\mu',\nu'})(x):= t^{\lambda}_{\mu,\nu}(x') t^{\lambda'}_{\mu',\nu'}(x''),
\end{eqnarray*}
see \cite{KS}. Furthermore, ($\mathbb{C}[G_q], \triangle_q, S_q, \epsilon_q$) is a dual Hopf algebra for $U_q(\g)$ equipped with
\begin{eqnarray*}
\triangle_q(t^{\lambda}_{\mu, \nu}) := \sum_{\kappa} t^{\lambda}_{\mu, \kappa} \otimes t^{\lambda}_{\kappa, \nu},\5 t^{\lambda}_{\mu, \nu}(S_q(x)) := (S_q(t^{\lambda}_{\mu, \nu}))(x),\5 \epsilon_q(t^{\lambda}_{\mu, \nu}) := t^{\lambda}_{\mu, \nu}(1),
\end{eqnarray*}
for all $x \in U_q(\g)$. Thus, we use the same symbols for the antipode and counit for $U_q(\g)$ and $\mathbb{C}[G_q]$. Using the identity $1 = \sum_{\mu} |\lambda, \mu \rangle \langle \lambda, \mu |$ we can write the left action in the form
\begin{eqnarray*}
\partial(x) t_{\mu \nu}^{\lambda} = ((t^{\lambda}_{\mu \nu})''(x)) (t^{\lambda}_{\mu \nu})'.
\end{eqnarray*}

\textbf{1.3.}\ The Haar state $h: \mathbb{C}[G_q] \rightarrow \mathbb{C}$ of the algebra $\mathbb{C}[G_q]$ is a positive $(h(t^* t) > 0)$ and left invariant functional, i.e. for all $t \in \mathbb{C}[G_q]$
\begin{eqnarray*}
(h \otimes \mathrm{id}) \triangle_q(t) = h(t).
\end{eqnarray*}
We shall assume the normalization $h(1) = 1$. The Haar state is faithful and therefore we can set an inner product   
\begin{eqnarray*}
\langle t,s \rangle = h(t^* s)
\end{eqnarray*} 
The Hilbert space $L^2(G_q)$ is the completion of $\mathbb{C}[G_q]$ with respect to the inner product. The GNS construction defines a faithful $*$-representation $\mathbb{C}[G_q] \rightarrow B(L^2(G_q))$. 

\section{Dirac Operators} 

We have a construction of a Hilbert space $L^2(G_q)$ and a faithful representation of the regular functions $\mathbb{C}[G_q]$ on $L^2(G_q)$. This system is equipped with a symmetry algebra $U_q(\g)$ which has a representation by unbouded densely defined operators on $L^2(G_q)$. In the following we describe two general methods how to equip this theory with an equivariant self adjoint operator which are $q$-deformations of Dirac operators. We give a  more detailed treatment for $SU_q(2)$, $U_q(2)$ and $SU_q(3)$. \\

\textbf{2.1.} We begin with the Lie theory. Consider the Hilbert space $L^2(G)$ and let $G$ act on it from the left. Then $L^2(G)$ decomposes into irreducible components so that by Peter-Weyl theorem $\mathbb{C}[G] = \bigoplus_{\lambda \in P_+} V^*_{\lambda} \otimes V_{\lambda}$ is a dense subset in $L^2(G)$. 

Let $\mathrm{cl}(\g)$ be the Clifford algebra affiliated with the vector space $\g$ and the Killing form. Let $\gamma: \g \rightarrow \mathrm{cl}(\g)$ denote the canonical embedding satisfying $\gamma(x)^2 = (x,x)\textbf{1}$. There exists a Lie algebra homomorphism $\widetilde{\mathrm{ad}}: \g \rightarrow \mathrm{cl}(\g)$ so that 
\begin{eqnarray}\label{cov}
\gamma([x,y]) = [\widetilde{\mathrm{ad}}(x), \gamma(y)]
\end{eqnarray}
which is given by
\begin{eqnarray*}
x \mapsto \widetilde{\mathrm{ad}}(x):=  \frac{1}{4} \sum_k \gamma(x_k) \gamma([x,x_k]) \in \mathrm{cl}(\g).
\end{eqnarray*}
Fix an irreducible representation $(\Sigma, s)$ for $\mathrm{cl}(\g)$. $(\Sigma, s)$ is called a spinor module. 

Fix a basis $\{x_i\}$ for the Lie algebra $\g$ so that $(x_i, x_j) = - \delta_{ij}$. The classical Dirac operator is defined by
\begin{eqnarray}\label{CD}
\mathcal{D} = \sum_k (  x_k \otimes \gamma(x_k) + \frac{1}{2} \otimes \gamma(x_k)\widetilde{\mathrm{ad}}(x_k) ) \in U(\g) \otimes \mathrm{cl}(\g)
\end{eqnarray}
The Hilbert space of square integrable sections of the spin bundle can be identified with $L^2(G) \otimes \Sigma$. The left action $\partial$ of $U(\g)$ on $L^2(G)$ makes $D:= (\partial \otimes s)\mathcal{D}$ an unbounded self adjoint operator on $L^2(G) \otimes \Sigma$. The Dirac operator is equivariant because it commutes with the image of the Lie algebra homomorphism $\g \rightarrow \g \otimes \mathrm{cl}(\g)$ defined by $x \mapsto x' \otimes \widetilde{\mathrm{ad}}(x'')$. This is the case because the generators of $\g$ and $\mathrm{cl}(\g)$ are covariant under the adjoint action of $\g$. \\

\textbf{2.2. Geometric Dirac Operator.} This operator is defined in \cite{NT} and its properties has been studied in \cite{NT2, NT3}. The term geometric here refers to the fact that the operator is suitable to define a geometric model for any quantum group deformation of a simple, simply connected and compact Lie group. We shall return to the geometric consideration in chapter 3. 

Let $\mathcal{D}$ be a classical Dirac operator on $G$. For each $\lambda \in P_+$ fix a representation of $U(\g)$ and $U_q(\g)$, a $*$-isomorphism $\phi: W^*(G_q) \rightarrow W^*(G)$ together with a unitary twist $F$ compatible with the isomorphism. Define the geometric Dirac operator on $G_q$ by 
\begin{eqnarray*}
\mathcal{D}_q = (\phi^{-1} \otimes \mathrm{id}) \Big( ( \mathrm{id} \otimes \widetilde{\mathrm{ad}})(F)\mathcal{D} ( \mathrm{id} \otimes \widetilde{\mathrm{ad}})(F^*) \Big) \in  U(G_q) \otimes \mathrm{cl}(\g) .
\end{eqnarray*}
$\mathcal{D}$ acts on the Hilbert space $L^2(G_q) \otimes \Sigma$ by 
\begin{eqnarray*}
D_q = (\partial \otimes s)\mathcal{D}_q.
\end{eqnarray*}
This operator is unitarily equivalent to the classical Dirac operator and therefore it has the same spectrum. $D_q$ is an unbounded and selfadjoint operator on $L^2(G_q) \otimes \Sigma$. Let $x \in W^*(G_q)$ act on $L^2(G_q) \otimes \Sigma$ by $x \mapsto \partial(x') \otimes s(\widetilde{\mathrm{ad}} \circ \phi (x))$. Using the corresponding property of the clasical Dirac operator it is straightforward to check that $D_q$ commutes with this action. \\

\textbf{2.3. Algebraic Dirac Operator.} This approach presented in \cite{H} is based on creating algebras $\mathfrak{L}_q(\g) \subset U_q(\g)$ and $\mathrm{cl}_q(\g) \subset B(\Sigma)$ which transform covariantly under the adjoint action of the Hopf algebra $U_q(\g)$. Then we can fix a suitable element in $U_q(\g) \otimes \mathrm{cl}_q(\g)$ which spans a singlet for the adjoint action and as a Hilbert space operator commutes with the representation of the symmetry algebra $U_q(\g)$ on $L^2(G_q) \otimes \Sigma$.       

Consider the adjoint representation $(V_{q,\rho}, \pi_{q,\rho})$ of $U_q(\g)$. Fix a basis $\{|\rho,n  \rangle: n \in I\}$ of $V_{q,\rho}$ and an orthonormal dual basis  $\{\langle \rho, m|: m \in I\}$ of $V^*_{q,\rho}$. The vector $\Omega = \sum_{n} |\rho,n  \rangle \otimes \langle \rho, n|$ spans the singlet of $V_{q,\rho} \otimes V^*_{q,\rho}$, i.e.
\begin{eqnarray*}
(\pi_{q,\rho} \otimes \pi^*_{q,\rho})(\triangle_q(x)) \Omega = \epsilon_q(x) \Omega, 
\end{eqnarray*}
for all $x \in U_q(\g)$. We would like to map $\Omega$ onto $\mathfrak{L}_q(\g) \otimes \mathrm{cl}_q(\g) \subset U_q(\g) \otimes B(\Sigma)$ by a module isomorphism and define the Dirac operator in its image.

We first define the covariant algebra $\mathrm{cl}_q(\g)$. In order to do this we explain the structure of the irreducible representations of $\mathrm{cl}(\g)$. Consider the Clifford algebra $\mathrm{cl}(\g)$ and let $(\Sigma,s)$ be a spinor module and  $s(\widetilde{\mathrm{ad}})$ the representation of $\g$ on $\Sigma$. The image of the embedding $\gamma: \g \rightarrow \mathrm{cl}(\g)$ forms an adjoint $\g$-module by \eqref{cov}. Denote by $\Psi$ the vector space $s(\gamma(\g)) \subset B(\Sigma) $. $\Psi$ spans an adjoint $\g$-module under the action 
\begin{eqnarray}\label{sady}
x \ad Y = [s(\widetilde{\mathrm{ad}}(x)), Y ], 
\end{eqnarray}
for $x \in \g$, $Y \in \Psi$. We equip the tensor product $B(\Sigma) \otimes B(\Sigma)$ with a structure of $\g$ module by applying the coproduct with the action \eqref{sady}. The multiplication $m: B(\Sigma) \otimes B(\Sigma) \rightarrow B(\Sigma)$ is a module homomorphism. If $(\Psi \otimes \Psi)_+ \subset \Psi \otimes \Psi$ denotes the reducible submodule spanned by the symmetric tensor products, then $m$ defines a homomorphism $(\Psi \otimes \Psi)_+ \rightarrow \mathbb{C}$ because $\Psi$ is spanned by $\mathrm{cl}(\g)$ representation matrices. As an algebra $\Psi$ generates $B(\Sigma)$. In terms of representation theory this has the following explanation. The module homomorphism $m$ can be applied in $\Psi$ to create new submodules of $B(\Sigma)$. Then the vector space $\Psi \oplus m(\Psi \otimes \Psi)$ spans a reducible subspace of $B(\Sigma)$. If this space is not isomorphic to $B(\Sigma)$ we need to apply $m$ again in it to create new submodules. After finite steps the morphism $m$ will create each irreducible component of $B(\Sigma)$. 

Recall that for $q = e^{i \pi h}$ with $h \in \mathbb{C} - \mathbb{Q}^*$ the category of finite dimensional $U_q(\g)$-modules, $\mathfrak{C}(\g,q)$, is equipped with a monoidal operation making it a braided monoidal category. The monoidal operation is the tensor product of modules and braiding can be defined using the $R$-matrix. The Drinfeld category, $\mathfrak{D}(\g,q)$, is a category of finite dimensional $\g$-modules. The standard tensor product with associativity morphism defined by the Drinfeld's associator $\Phi_{\mathrm{KZ}}$ defines a braided monoidal structure. There is a $\mathbb{C}$-linear braided  monoidal equivalence between the catories $\mathfrak{C}(\g,q)$ and $\mathfrak{D}(\g,q)$, \cite{KL1},\cite{KL2},\cite{NT0}. The modules with the same highest weights in $\mathfrak{C}(\g,q)$ and $\mathfrak{D}(\g,q)$ correspond to each other under the equivalence. The categoties $\mathfrak{C}(\g,q)$ and $\mathfrak{D}(\g,q)$ are also rigid. Especially the dual modules provide the dual objects.  

In a rigid monoidal category there is a morphism $ \mathrm{ev}_{X} : X^* \otimes X \rightarrow \textbf{1}$ for each object $X$. In the categories $\mathfrak{C}(\g,q)$ and $\mathfrak{D}(\g,q)$ these morphisms coincide, up to a multiplicative constant, with the standard basis independent pairing $X^* \otimes X \rightarrow \textbf{1}$ of the category of vector spaces \cite{NT0}. Let $F: \mathfrak{D}(\g,q) \rightarrow \mathfrak{C}(\g,q)$ denote the category equivalence which is also a functor of rigid monoidal categories.  Then we have a morphism 
$ F(\mathrm{ev}_{X}) : F(X^*) \otimes F(X) \rightarrow \textbf{1}$. There is a unique morphism $\alpha: F(X^*) \rightarrow F(X)^*$ so that 
\begin{eqnarray*}
F(\mathrm{ev}_{X}) = \mathrm{ev}_{F(X)} \circ (\alpha \otimes \mathrm{id}) : F(X^*) \otimes F(X) \rightarrow \textbf{1}.
\end{eqnarray*}
Since $F$ is a functor between rigid monoidal categories $\alpha$ is an isomoprhism \cite{DM} (Proposition 1.9.). The modules $F(X^*)$ and $F(X)^*$ are isomorphic. Therefore $\mathrm{ev}_{F(X)}$ can be considered as a map $F(X^*) \otimes F(X) \rightarrow \textbf{1}$ and then $F(\mathrm{ev}_{X})$ and $\mathrm{ev}_{F(X)}$ coincide, up to a multiplicative constant, in $\mathfrak{C}(\g,q)$. If $(X, \pi_q)$ is in $\mathfrak{C}(\g,q)$ we consider $B(X)$ as a module with the usual action
\begin{eqnarray*}
x \ad_q Y := \pi_q(x') Y \pi_q(S_q(x'')) 
\end{eqnarray*}
for all $x \in U_q(\g)$, $Y \in B(X)$. The matrix multiplication $m: B(X) \otimes B(X) \rightarrow B(X)$ defines a module homomorphism in $\mathfrak{C}(\g,q)$ and $\mathfrak{D}(\g,q)$ and the equivalence of monoidal categories implies that its domain $B(X) \otimes B(X)$ has the same decomposition into irreducible components in both categories. Since $B(X) \simeq X \otimes X^*$ we can write $m = \mathrm{id} \otimes \mathrm{ev}_{X} \otimes \mathrm{id}$. Then $F(m) = \mathrm{id} \otimes F(\mathrm{ev}_{X}) \otimes \mathrm{id}$ coincides with the multiplication morphism $m$ in $\mathfrak{C}(\g,q)$ up to a constant. From the equivalence of categories we know that the irreducible components of $m(B(X) \otimes B(X))$ in $\mathfrak{D}(\g,q)$ and $ F(m)(F(B(X)) \otimes F(B(X)))$ in $\mathfrak{C}(\g,q)$ are the same. The objects $F(B(X))$ can be identified with $B(X)$  but it is understood to carry a representation of $U_q(\g)$. We conclude that the domain $B(X) \otimes B(X)$ and the image $m(B(X) \otimes B(X))$ have the same decomposition into irreducible components in both categories $\mathfrak{D}(\g,q)$ and $\mathfrak{C}(\g,q)$. 

Let $\Sigma$ in $\mathfrak{C}(\g,q)$ correspond to the spinor module under the category equivalence. From the module $B(\Sigma)$ in $\mathfrak{C}(\g,q)$ we can pick a component $\Psi_q$ which corresponds to the component $\Psi$ in $\mathfrak{D}(\g,q)$ and $m: (\Psi_q \otimes \Psi_q)_+ \rightarrow \textbf{1}$ is onto where $(\Psi_q \otimes \Psi_q)_+ \subset \Psi_q \otimes \Psi_q$ has the same isotypic components as $(\Psi \otimes \Psi)_+$ has in $\mathfrak{D}(\g,q)$. We denote by $\mathrm{cl}_q(\g)$ the associative algebra generated by $\Psi_q$. Again the multiplication morphism can be applied to generate the endomorphism algebra $B(\Sigma)$ from the subspace $\Psi_q$ because the image of $m$ has the same components as in $\mathfrak{D}(\g,q)$. Therefore $\mathrm{cl}_q(\g)$ is the algebra $B(\Sigma)$. We shall demonstrate how to find the covariant generators for $\mathrm{cl}_q(\g)$ in the cases $\g = \mathfrak{su_2}$ and $\g = \mathfrak{su_3}$ below. Since $V_{q,\rho}^* \simeq V_{q,\rho}$ we can fix a module isomorphism $\sigma: V_{q,\rho}^* \rightarrow \Psi_q$.

Next we recall that there exists a subspace, $\mathfrak{L}_q(\g) \subset U_q(\g)$, called a quantum Lie algebra which transforms covariantly under the adjoint action of $U_q(\g)$ on itself 
\begin{eqnarray*}
x \bad y = x''y S(x'). 
\end{eqnarray*}
This action is actually an opposite adjoint action, however, it is exactly what we need to make the Dirac operator equivariant, see \cite{H}. Furthermore $\mathfrak{L}_q(\g)$ reduces to $\g$ in the classical limit $q \rightarrow 1$. For construction see \cite{De}. Denote by $\theta: V_{q, \rho} \rightarrow \mathfrak{L}_q(\g)$ a module isomorphism.

The Dirac operator is defined by 
\begin{eqnarray*}
\mathcal{D}'_q = (\theta \otimes \sigma)\Omega \in  U_q(\g) \otimes \mathrm{cl}_q(\g),
\end{eqnarray*}
This operator commutes with the image of the homomorphism $x \mapsto (\mathrm{id} \otimes \pi_q )\triangle_q(x)$, see \cite{H}. We can define a Hilbert space operator by setting
\begin{eqnarray*}
\mathfrak{D}_q = (\partial \otimes \mathrm{id}) \mathcal{D}'_q
\end{eqnarray*}
which is an equivariant self adjoint operator on $L^2(G_q) \otimes \Sigma$. \\

\textbf{2.4. Dirac Operators on $SU_q(2)$.} Let us choose the generators $\{j_{\pm}, j_0\}$ of $\mathfrak{su_2}$ so that
\begin{eqnarray*}
[j_0, j_{\pm} ] = \pm j_{\pm}, \5 [j_+, j_-] = 2j_0, \5 C = j_+ j_- + j_0(j_0 + 1) := j(j+1)
\end{eqnarray*}
where the Casimir operator $C$ is also given. The irreducible representations $\{(V_l, \pi_l): l \in \frac{1}{2} \mathbb{N}_0 \}$ are given by 
\begin{eqnarray*}
\pi_l (j_{\pm}) |l,m \rangle = \sqrt{l(l+1)-m(m \pm 1)} |l,m \pm 1 \rangle ,\5  \pi_l (j_0) |l,m \rangle = m |l,m \rangle. 
\end{eqnarray*}
where the basis is chosen by $\{ |l,m \rangle: -l \leq m \leq l\}$ for each $V_l$. The Killing form is normalized so that the vectors 
\begin{eqnarray*}
x_1 = j_+ + j_-, \5  x_2 = -i(j_+ - j_-), \5 x_3 = 2 j_0
\end{eqnarray*}
form an orthonormal basis of $\g$. The representations of the algebras $\mathrm{cl}(\mathfrak{su_2})$ and $\mathfrak{su_2}$ on $\Sigma = V_{\frac{1}{2}}$ are
\begin{eqnarray*}
s: \gamma(x_i) \mapsto  \pi_{\frac{1}{2}}(x_i), \5 \widetilde{\mathrm{ad}}(x_i) \mapsto \pi_{\frac{1}{2}}(x_i).
\end{eqnarray*}
The classical Dirac operator \eqref{CD} corresponding to these choices is defined by
\begin{eqnarray}\label{const}
D = (\partial \otimes s)\mathcal{D} = 2 \partial \begin{pmatrix} 
  j_0 & j_- \\
  j_+ & -j_0 \end{pmatrix} + \frac{3}{2} \mathbf{1}.
\end{eqnarray}
$D$ acts on $L^2(SU(2)) \otimes V_{\frac{1}{2}}$. 

Fix the unitary representations $(V_{q,l}, \pi_{q,l})$ of $U_q(\mathfrak{su_2})$ by 
\begin{eqnarray*}
\pi_{q,l}(k)|l,m \rangle &=& q^m |l,m \rangle  \\ \nonumber
\pi_{q,l}(e)|l,m \rangle &=& \sqrt{[l-m][l+m+1]} |l,m+1 \rangle \\ \nonumber 
\pi_{q,l}(f)|l,m \rangle &=& \sqrt{[l-m+1][l+m]} |l,m-1 \rangle 
\end{eqnarray*}
for all $l \in \frac{1}{2} \mathbb{N}_0$ where $[n] := (q^n - q^{-n})(q-q^{-1})^{-1}$. Choose an algebra $*$-isomorphism $\phi_l: B(V_{q,\lambda}) \rightarrow B(V_{\lambda})$ for each $l \in P_+$ by  (cf. \cite{Cu}) 
\begin{eqnarray*}
& &\phi_l(\pi_{q,l}(e)) =  \sqrt{\frac{[\pi_l(j)-\pi_l(j_0)+1][\pi_l(j)+\pi_l(j_0)]}{\pi_l(j)(\pi_l(j)+1)-\pi_l(j_0)(\pi_l(j_0) - 1)}} \pi_l(j_{+}) \\
& &\phi_l(\pi_{q,l}(f)) = \sqrt{\frac{[\pi_l(j)-\pi_l(j_0)][\pi_l(j)+\pi_l(j_0)+1]}{\pi_l(j)(\pi_l(j)+1)-\pi_l(j_0)(\pi_l(j_0) + 1)}} \pi_l(j_-), \5 \phi_l(\pi_{q,l}(k)) = q^{\pi_l(j_0)}.
\end{eqnarray*}
These define a $*$-isomorphism $\phi: W^*(G_q) \rightarrow W^*(G)$ and we fix a unitary twist $F$ compatible with $\phi$. We have $\pi_{l} \circ \phi = \pi_{q,l}$ for each $l$. 

For now we set the constant operator $(3/2)\textbf{1}$ in \eqref{const} to zero and denote by $\tilde{D}_q$ the geometric Dirac operator without this constant term: $D_q = \tilde{D}_q + (3/2)\textbf{1}$. In the present example we have $s(\gamma(x_i)) = s((\widetilde{\mathrm{ad}})(x_i)) = \pi_{\frac{1}{2}}(x_i)$ for all $x_i \in \g$. It was noted in \cite{NT} that in this case we can use the relation 
\begin{eqnarray*}
(\phi \otimes \phi)(R^* R) = F q^{T}F^*, \5 T:= \sum_k x_k \otimes x_k
\end{eqnarray*}
to write  
\begin{eqnarray*}
q^{\tilde{D}_q} &=& (\partial \circ \phi^{-1} \otimes \pi_{\frac{1}{2}})(F q^{T}F^*) = (\partial \otimes \pi_{q,\frac{1}{2}})(R^* R) \\
 &=& \partial \Big[\begin{pmatrix} 
  k^2 & 0 \\
  0 & -k^2 \end{pmatrix} + (q-q^{-1}) \begin{pmatrix} 
  (1-q^{-2})fe &\ \ q^{-\frac{1}{2}}fk^{-1} \\
  q^{-\frac{1}{2}}k^{-1}f &\ \ 0 \end{pmatrix} \Big].
\end{eqnarray*}
An explicit formula for the $R$-matrix used in the calculation can be found from \cite{KR}. The relation $s \circ \gamma = s \circ \widetilde{\mathrm{ad}}$ does not hold for a general $\g$ and more advanced methods are needed in order to find an explicite formula. 

The algebraic operator on $SU_q(2)$ was constructed in \cite{H}. Now $(\Sigma,\pi_{q})$ is $(V_{q,\frac{1}{2}},\pi_{q,\frac{1}{2}})$. The adjoint representation of $U_q(\mathfrak{su_2})$ is $(V_{q,1}, \pi_{q,1})$. As a module we have 
\begin{eqnarray*}
B(\Sigma) = V_{q,1} \otimes V^*_{q,1} \simeq V_{q,1} \otimes V_{q,1} \simeq V_{q,0} \oplus V_{q,1} \oplus V_{q,3}.
\end{eqnarray*}
Therefore the adjoint representation $V_{q,1} \subset B(\Sigma)$ does not have multiplicities and its generators span the covariant module algebra $\Psi_q$ which generates the algebra $\mathrm{cl}_q(\mathfrak{su_2})$. The highest weight condition gives immediately the following basis for $\Psi_q$ which is unique up to scaling
\begin{eqnarray*}
s(\psi_1) = \begin{pmatrix} 
  0 & \sqrt{q} \\
  0 & 0 \end{pmatrix}, \5 s(\psi_0) = - \frac{1}{\sqrt{[2]}}\begin{pmatrix} 
  q^{-1} & 0 \\
  0 & -q \end{pmatrix}, \5 s(\psi_{-1}) = \begin{pmatrix} 
  0 & 0 \\
  -\sqrt{q^{-1}} & 0 \end{pmatrix}.
\end{eqnarray*}
It is straightforward to check that $m: \Psi_{q,+} = V_3 \oplus V_0 \rightarrow \mathbb{C}$ holds and the $q$-deformed Clifford relations are
\begin{eqnarray*}
& &\psi_{q,1} \psi_{q,1} = \psi_{q,-1} \psi_{q,-1} = 0 \\
& & q^{-1} \psi_{q,1} \psi_{q,0} + q \psi_{q,0} \psi_{q,1} = 0 \\
& &q^{-2} \psi_{q,1} \psi_{q,-1} + [2] \psi_{q,0} \psi_{q,0} + q^2 \psi_{q,-1} \psi_{q,1} = 0 \\
& &\psi_{q,0} \psi_{q,-1} + q^2 \psi_{q,-1} \psi_{q,0} = 0 \\
& &\psi_{q,1} \psi_{q,-1} + \psi_{q,-1} \psi_{q,1} = -1,
\end{eqnarray*}
where $\psi_{q,i} = |1,i \rangle \in B(\Sigma)$. 

The isomorphism $V \rightarrow \mathfrak{L}(\mathfrak{su_2})$ is defined by 
\begin{eqnarray*}
\theta(|1,1 \rangle) = k^{-1}e, \5 \theta(|1,0 \rangle) = \frac{1}{\sqrt{[2]}}(q^{-1}fe - qef),\5 \theta(|1,-1 \rangle) = - k^{-1}f.
\end{eqnarray*}
Therefore we find the following  Dirac operator 
\begin{eqnarray*}
\mathfrak{D}_q = \partial \begin{pmatrix} 
  ef-q^{-2}fe\  & \  q^{-\frac{1}{2}} [2] k^{-1}f \\
  q^{\frac{1}{2}}[2] k^{-1}e\ &\ -q^2 ef + fe \end{pmatrix}.
\end{eqnarray*}

There is a fundamental relationship between the geometric and algebraic approach
\begin{eqnarray*}
\mathfrak{D}_q = [\tilde{D}_q] = \frac{q^{\tilde{D}_q} - q^{-\tilde{D}_q}}{q-q^{-1}}
\end{eqnarray*}
which can be checked by using the formula $q^{-\tilde{D}_q} = (\partial \otimes \pi_{q,\frac{1}{2}})(R^{-1} (R^*)^{-1})$. This explains the very different spectral behaviours of these operators. \\

\textbf{2.5. Geometric Dirac Operator on $U_q(2)$.} The Lie algebra $\mathfrak{u_2}$ is spanned by $x_i$ ($ 0 \leq i \leq 3$) so that $x_0$ is central and $x_1, x_2, x_3$, defined as above, span the subalgebra $\mathfrak{su_2}$. Let us fix the normalization of the nondegenerate bilinear form so that these $x_i$ form an orthonormal basis. The irreducible finite dimensional representations are parametrized by the pairs $(l,c)$ where $l$ is the highest weight of the subalgebra $\mathfrak{su_2}$ and $c$ fixes the action of the center. Furthermore, $l$ and $c$ are both integers or both half integers. Denote by $P_+$ the set of such pairs. 

The $q$-deformed algebra $U_q(\mathfrak{u_2})$ is defined by adding the linearly independent generators $\xi_1, \xi_1^{-1}, \xi_2$ and $\xi_2^{-1}$ to the algebra $U_q(\mathfrak{su_2})$ so that $k = \xi_1 \xi_2^{-1}$ and the element $\xi_1 \xi_2$ is central. The extension of the Hopf structure is defined by
\begin{eqnarray*}
\triangle_q(\xi_i) = \xi_i \otimes \xi_i,\5  S_q(\xi_i) = \xi_i^{-1}, \5 \epsilon_q(\xi_i) = 1,\5 i=1,2.
\end{eqnarray*}
Since $U_q(\mathfrak{u_2})$ differs from $U_q(\mathfrak{su_2})$ only by an element in the center, the twist $F$ is defined as in the case  $U_q(\mathfrak{su_2})$. Again the highest weight modules are parametrized by the pairs $(l,c) \in P_+$ (both integers or half integers) where $l$ is the highest weight of the subalgebra $U_q(\mathfrak{su}_2)$ and the central element $\xi_1 \xi_2$ acts by $q^c$. 

The Clifford algebra $\mathrm{cl}(\mathfrak{u_2})$ has a four dimensional irreducible representation $\hat{\Sigma}$ given by
\begin{eqnarray*}
s(\gamma(x_0)) = \begin{pmatrix} 
  0 & \textbf{1} \\
  \textbf{1} & 0 \end{pmatrix}, \5 s(\gamma(x_n)) = i  \begin{pmatrix} 
  0 & \pi_{\frac{1}{2}}(x_n) \\
  -\pi_{\frac{1}{2}}(x_n) & 0 \end{pmatrix}, \5 1 \leq n \leq 3. 
\end{eqnarray*}
The corresponding representation $s(\widetilde{\mathrm{ad}})$ of $\mathfrak{u_2}$ has two irreducible components $\hat{\Sigma} = V^+_{(\frac{1}{2},0)} \oplus V^-_{(\frac{1}{2},0)}$
\begin{eqnarray*}
s(\widetilde{\mathrm{ad}}(x_0))= 0, \5 s(\widetilde{\mathrm{ad}}(j_x)) =   \begin{pmatrix} 
  \pi_{\frac{1}{2}}(j_x) & 0 \\
  0 & \pi_{\frac{1}{2}}(j_x) \end{pmatrix}, \5 x \in \{\pm, 0\},
\end{eqnarray*}
where we have fixed the action of the center to be zero. Therefore, if we denote by $D$ and $D_q$ the Dirac operators on $SU(2)$ and $SU_q(2)$ we get 
\begin{eqnarray*}
\hat{D} =  \begin{pmatrix} 
  0 & iD +  \partial(x_0) \\
  -iD +  \partial(x_0) & 0  \end{pmatrix},\5 \hat{D}_q =  \begin{pmatrix} 
  0 & iD_q + \partial(x_0) \\
  -iD_q + \partial(x_0) & 0  \end{pmatrix}.
\end{eqnarray*}
We return to this construction in Chapter 3.\\

\textbf{2.6. Algebraic Dirac Operator on $SU_q(3)$.} An irreducible representation $(\Sigma,s)$ of $\mathrm{cl}(\mathfrak{su}_3)$ is $16$-dimensional and $(\Sigma, s(\widetilde{\mathrm{ad}}))$ splits into two components, $\Sigma_+$ and $\Sigma_-$, both isomorphic to the adjoint representation $(V_{\rho}, \pi_{\rho})$. The basis of $\Sigma:= \Sigma_+ \oplus \Sigma_-$ can be chosen so that the representations are of the form
\begin{eqnarray*}
 s(\widetilde{\mathrm{ad}})(x)  = \begin{pmatrix} 
  \pi_{\rho}(x) & 0  \\
  0 & \pi_{\rho}(x)  \end{pmatrix},\5  s(\gamma(x_i)) = \begin{pmatrix} 
  0 & \psi^+_i  \\
  \psi^-_i & 0  \end{pmatrix}
\end{eqnarray*}
for all $x \in \mathfrak{su}_3$ and for all generators $\gamma(x_i)$ of $\mathrm{cl}(\mathfrak{su}_3)$. The off diagonal operators $\psi^{\pm}_i:  \Sigma_{\mp} \rightarrow \Sigma_{\pm}$ transform covariantly as an adjoint representation under the action $(x, \psi^{\pm}_i) \mapsto [\pi_{\rho}(x), \psi^{\pm}_i]$.

The deformed representations preserve the structural zeros and therefore we can assume that the representation of $U_q(\mathfrak{su_3})$ and the covariant space $\Psi_q$ are operators on $\Sigma$ of the form 
\begin{eqnarray*}
 \pi_q(x)  = \begin{pmatrix} 
  \pi_{q,\rho}(x) & 0  \\
  0 & \pi_{q,\rho}(x)  \end{pmatrix},\5  \psi_{q,i} = \begin{pmatrix} 
  0 & \psi^+_{q,i}  \\
  \psi^-_{q,i} & 0  \end{pmatrix}
\end{eqnarray*}
We use the Gelfand-Tsetlin basis for the representation $\pi_{q,\rho}$, i.e.,
\begin{eqnarray*}
& &\pi_{q,\rho}(e_1) = e_{12} + \sqrt{[2]}e_{35} + \sqrt{[2]}e_{56} + e_{78}, \\
& &\pi_{q,\rho}(e_2) = e_{13} + \sqrt{[3]/[2]}e_{24} + (\sqrt{[2]})^{-1} e_{25} + \sqrt{[3]/[2]}e_{47} + (\sqrt{[2]})^{-1} e_{57} + e_{68},\\
& &\pi_{q,\rho}(k_1) = q^{\frac{1}{2}}e_{11} + q^{-\frac{1}{2}}e_{22} + qe_{33} + e_{44} + e_{55} + q^{-1}e_{66} + q^{\frac{1}{2}}e_{77} + q^{-\frac{1}{2}}e_{88}, \\
& &\pi_{q,\rho}(k_2) = q^{\frac{1}{2}}e_{11} + qe_{22} + q^{-\frac{1}{2}}e_{33} + e_{44} + e_{55} + q^{\frac{1}{2}}e_{66} + q^{-1}e_{77} + q^{-\frac{1}{2}}e_{88}, \\
& &\pi_{q,\rho}(f_1) = (\pi_{q,\rho}(e_1))^{\dagger}, \5 \pi_{q,\rho}(f_2) = (\pi_{q,\rho}(e_2))^{\dagger}.
\end{eqnarray*}
The adjoint action of $U_q(\mathfrak{su_3})$ on $\Psi_q \subset B(\Sigma)$ is given as 
\begin{eqnarray}\label{adact}
x \ad_q \psi_{q,i} = \begin{pmatrix} 
  0 & \pi_{q, \rho}(x) \ad_q \psi^+_{q,i}   \\
  \pi_{q, \rho}(x) \ad_q \psi^-_{q,i} & 0  \end{pmatrix}.
\end{eqnarray}
In order to make $\psi_{q,i}: 1 \leq i \leq 8$ covariant under \eqref{adact} we need to make the off diagonal operators $\psi^{\pm}_{q,i}$ covariant under the representation 
\begin{eqnarray*}
(x, \psi^{\pm}_{q,i}) \mapsto \pi_{q,\rho}(x) \ad_q \psi^{\pm}_{q,i} = \pi_{q,\rho}(x')\psi^{\pm}_{q,i}\pi_{q,\rho}(S_q(x'')).
\end{eqnarray*}
As a module, the spaces of off diagonal opearators are of the form $V_{q,\rho} \otimes V^*_{q,\rho} \simeq V_{q,\rho} \otimes V_{q,\rho}$. The tensor product has a reduction into irreducible components
\begin{eqnarray}\label{tensor}
V_{q,\rho} \otimes V_{q,\rho} \simeq V_{q,2 \rho} \oplus V_{q,0} \oplus  2V_{q,\rho} \oplus V_{q,\rho + \xi_1 - \xi_2} \oplus V_{q,\rho + \xi_2- \xi_3},
\end{eqnarray}
where we have written the simple roots by $\alpha_i = \xi_i - \xi_{i+1}: 1 \leq i \leq 2$. Since the adjoint representation occurs with a  multiplicity two in \eqref{tensor} this condition does not fix the operators even up to a multiplicative constant. A general form for the covariant operators is 
\begin{eqnarray*}
\psi_{q,1}^{\pm} &=& e_{14} + b_{\pm}e_{15} + q^{-\frac{1}{2}}\sqrt{[2]}b_{\pm} e_{26} + (\sqrt{q[2]})^{-1}(b_{\pm} + \sqrt{[3]})e_{37} \\
&+& (q[2])^{-1}(\sqrt{[3]}b_{\pm}-1)e_{48} + (q[2])^{-1}(\sqrt{[3]} + b_{\pm})e_{58} \\
\psi_{q,2}^{\pm} &=& -q^{-3/2}\sqrt{[2]}b_{\pm}e_{13} + e_{24} - q^{-2} b_{\pm} e_{25} - (q^2[2])^{-1}(\sqrt{[3]}b_{\pm}-1)e_{47} \\
&+& ([2])^{-1} (\sqrt{[3]} + b_{\pm}) e_{57} + (\sqrt{q[2]})^{-1} (\sqrt{[3]} + b_{\pm})e_{68} \\
\psi_{q,3}^{\pm} &=& -(\sqrt{q^3[2]})^{-1}(\sqrt{[3]} + b_{\pm}) e_{12} - (q[2])^{-1}(q^{-2} + \sqrt{[3]}b_{\pm})e_{34} + (q[2])^{-1}(q^2 b_{\pm} - \sqrt{[3]})e_{35} \\
&+& ([2])^{-1}(q^{-2} + \sqrt{[3]}b_{\pm})e_{46} + ([2])^{-1}(b_{\pm} - q^{-2}\sqrt{[3]}) e_{56} + q^{-\frac{1}{2}} \sqrt{[2]}b_{\pm} e_{78} \\
\psi_{q,4}^{\pm} &=& (q^3 [2])^{-1}(\sqrt{[3]}b_{\pm}- 1) e_{11}  + (q^3 [2])^{-1} (\sqrt{[3]}b_{\pm}- 1) e_{22} - (q^3[2])^{-1}  (1 + q^2\sqrt{[3]}b_{\pm})e_{33} \\
&+& (1- (q^3 [2])^{-1}(1 - \sqrt{[3]}b_{\pm}))e_{44} - (q^3[2])^{-1}  (1 + q^2\sqrt{[3]}b_{\pm})e_{55} \\
&-& (q^3[2])^{-1}(1 + q^2\sqrt{[3]}b_{\pm})e_{66} + e_{77} + e_{88} \\
\psi_{q,5}^{\pm} &=& (q^3[2])^{-1} (b_{\pm} + \sqrt{[3]})e_{11} - (q[2])^{-1} (b_{\pm} + \sqrt{[3]})e_{22} - (q^3[2])^{-1} (q^2b_{\pm}-\sqrt{[3]})e_{33} \\
&-& (q^3[2])^{-1} (q^2 \sqrt{[3]}b_{\pm} + 1)(e_{45} + e_{54}) + ([2])^{-1}(q-q^{-1})(b_{\pm} - q^{-2}\sqrt{[3]})e_{55} \\
&+& (q[2])^{-1}(q^2 b_{\pm} - \sqrt{[3]})e_{66} -q^{-2}b_{\pm}e_{77} + b_{\pm} e_{88} \\
\psi_{q,6}^{\pm} &=& (q^{5/2}\sqrt{[2]})^{-1}(\sqrt{[3]} + b_{\pm})e_{21} + (q^4[2])^{-1}(q^2 \sqrt{[3]}b_{\pm} + 1 )e_{43} + (q^2 [2])^{-1}(\sqrt{[3]}-q^2 b_{\pm})e_{53} \\
&-& (q^3[2])^{-1} (q^2 \sqrt{[3]}b_{\pm} + 1)e_{64} + (q^3[2])^{-1}( \sqrt{[3]} - q^2 b_{\pm})e_{65} - q^{-3/2} \sqrt{[2]}b_{\pm} e_{87} \\
\psi_{q,7}^{\pm} &=& q^{-5/2}\sqrt{[2]}b_{\pm} e_{31} -q^{-1} e_{42} + q^{-3} b_{\pm} e_{52} + (q^3 [2])^{-1}(\sqrt{[3]}b_{\pm} - 1)e_{74} \\
&-& (q [2])^{-1}(\sqrt{[3]} + b_{\pm})e_{75} - (q^{3/2} \sqrt{[2]})^{-1}(\sqrt{[3]} + b_{\pm})e_{86} \\
\psi_{q,8}^{\pm} &=& q^{-2} e_{41} + q^{-2}b_{\pm} e_{51} + q^{-5/2}\sqrt{[2]}b_{\pm}e_{62} + (q^{5/2}\sqrt{[2]})^{-1}(\sqrt{[3]} + b_{\pm})e_{73} \\
&+& (q^3[2])^{-1}(\sqrt{[3]}b_{\pm} - 1)e_{84} + (q^3[2])^{-1}(\sqrt{[3]} + b_{\pm})e_{85}
\end{eqnarray*}
where $b_{\pm}$ are free complex parameters. The operators $\psi^{\pm}_{q,i}$ are considered as $8 \times 8$ matrices which live in the off diagonal blocks of $16 \times 16$ matrices. We have fixed the scaling. 

Denote by $\Psi_q$ the vector space spanned by the operators $\psi_{q,i} : 1 \leq i \leq 8$. $\Psi_q$ is an adjoint module under the action \eqref{adact}. Also the tensor product $\Psi_q \otimes \Psi_q$ reduces according to \eqref{tensor}. The module isomorphic to $V_{q,2 \rho}$ in $\Psi_q \otimes \Psi_q$ reduces to a symmetric module in the classical case $q=1$. The vector $\psi_{q,1} \otimes \psi_{q,1}$ is the highest weight vector. The condition 
\begin{eqnarray*}
m: \psi_{q,1} \otimes \psi_{q,1} \mapsto \psi_{q,1}\psi_{q,1} = 0 
\end{eqnarray*}
leads to the constraint
\begin{eqnarray*}
b_- = \frac{1- \sqrt{[3]}b_+}{\sqrt{[3]} + b_+}. 
\end{eqnarray*}
Then the covariance implies that $\Psi_q \otimes \Psi_q \supset V_{q,2 \rho} \subset \mathrm{Ker}(m)$.

Another irreducible component of $2 V_{q,\rho} \subset \Psi_q \otimes \Psi_q$ reduces to a symmetric module for $q = 1$. Up to scaling a highest weight vector of an adjoint module in $ \Psi_q \otimes \Psi_q$ is of the form 
\begin{eqnarray*}
\omega_{\rho}(b_+,z) &=& (\sqrt{[3]})^{-1}(q^3[2]z - 1)\psi_{q,1} \otimes \psi_{q,4}+ (\sqrt{[3]})^{-1}(q^{-3}[2] - z)\psi_{q,4} \otimes \psi_{q,1} \\
&+& \psi_{q,1} \otimes \psi_{q,5} + z \psi_{q,5} \otimes \psi_{q,1} - q^{-3/2}\sqrt{[2]}  \psi_{q,2} \otimes \psi_{q,3} - q^{3/2}\sqrt{[2]}z \psi_{q,3} \otimes \psi_{q,2}.
\end{eqnarray*}
where $z \in \mathbb{C}$ is arbitrary. The condition $m: \omega_{\rho}(b_+,z) \mapsto 0$  leads to the following solutions for $z$ and $b_+$
\begin{eqnarray*}
z' = \frac{q^6 + q^2 +1}{(q^6 + q^4 + 1)q^4}, \5 b'_+ = \frac{q^2-1 \pm q [2] \sqrt{-[3]}}{2q^2 \sqrt{-[3]}}
\end{eqnarray*}
and 
\begin{eqnarray*}
z'' = \frac{q^6 + q^4 +1}{(q^6 + q^2 + 1)q^2}, \5 b''_+ = \frac{-q^2+1 \pm q [2] \sqrt{-[3]}}{2\sqrt{-[3]}}.
\end{eqnarray*}
In the classical case $q=1$ we could only have one possible value for $z$, namely $z=1$, which would give the highest weight vector for the symmetric component of $2V_{\rho}$. We would get two isomorphic irreducible $\mathrm{cl}(\mathfrak{su_3})$ modules. In the $q$-deformed case all the four possibilities reduce to these classical cases and therefore doubles the number of choices for the basis of the covariant generators. However, the algebras these generate are isomorphic. 

Finally we need to consider the trivial module $V_0 \subset \Psi_q \otimes \Psi_q$ which is the last submodule to study that reduces to a symmetric module for $q \rightarrow 1$. $V_0$ is spanned by 
\begin{eqnarray*}
\omega_0 &=& q^2\psi_{q,1} \otimes \psi_{q,8} + q^{-2} \psi_{q,8} \otimes \psi_{q,1} - q \psi_{q,2} \otimes \psi_{q,7} - q^{-1} \psi_{q,7} \otimes \psi_{q,2} \\
&-&  q \psi_{q,3} \otimes \psi_{q,6} - q^{-1} \psi_{q,6} \otimes \psi_{q,3} + \psi_{q,4} \otimes \psi_{q,4} + \psi_{q,5} \otimes \psi_{q,5}
\end{eqnarray*}
We have 
\begin{eqnarray*} 
m(\omega_0) = \frac{[4](2b_+ +(1-b_+^2)\sqrt{[3]})}{q^3(b_+ + \sqrt{[3]})} \textbf{1}.
\end{eqnarray*}
Therefore we have found the covariant generators $\Psi_q$ and the algebra $\mathrm{cl}_q(\mathfrak{su_3})$. To write down a formula for the Dirac operator one needs a quantum Lie algebra. An explicit construction can be found in \cite{De}. 

\section{Spectral Geometry} 

\textbf{3.1.} A unital spectral triple $(A, D, \h)$ consists of the following pieces of data: a unital associative $*$-algebra $A$ with a faithful $*$-representation $\rho$ on a separable Hilbert space $\h$. The operator $D$ is an unbounded self-adjoint operator with a dense domain in $\h$ such that $[D, \rho(t)]$ defined in the domain of $D$ extends to an bounded operator on $\h$ for all $t \in A$. The dimension of a spectral triple is the smallest integer so that $(1+D^2)^{-\frac{n}{2}} \in L_{1+}(\h)$ (The first Dixmier ideal). In the even dimensional case, there exists a chirality operator $\gamma \in B(\h)$ satisfying: $\gamma D + D \gamma = 0$ and $[\gamma, \rho(t)] = 0$ for all $t \in A$.

In the quantum group model under consideration we choose $A = \mathbb{C}[G_q]$. We use the Haar state to complete $A \otimes \Sigma$ to a Hilbert space $\h$ and the GNS construction gives us a faithtul $*$-representation of $\mathbb{C}[G_q]$ on $\h$. There exists a natural action of both Dirac operators $D_q$ and $\mathfrak{D}_q$ on $\h$. We recall the following theorem from \cite{NT}. \\

\noindent \textbf{Theorem.} Let $G$ denote a simple, simply connected and compact Lie group. Then $(\mathbb{C}[G_q], D_q, \h)$ is a spectral triple whose dimension matches with the dimension of the Lie group $G$. There exists an explicit formula for the chirality $\gamma$. \\

Smoothness of a noncommutative space is described by a property called regularity. Denote $\delta: \rho(t) \mapsto [|D|,  \rho(t)]$. A spectral triple is regular if the algebra generated by $\rho(t), [D, \rho(t)]$ for $t \in A$ is in the domain of $\delta^k$ for each $k \geq 0$. It is not known if the geometric spectral triples in the above theorem satisfy regularity. Besides the axioms given here there are suplementary axioms for the full description of NC Riemannian spin geometry \cite{Con}. However, in the quantum group case all these cannot be assumed to be fullfilled in their original form, see the discussion of \cite{DV}.

We first study the spectral triple of $SU_q(2)$ with details. It turns out that the isospectral deformation \cite{DV} is, up to one convention, the spectral triple of \cite{NT} associated to the Dirac operator we discussed in chapter $2$. The Fredholm modules associated to the Algebraic and geometric approaches turn out to be homotopic. We use the details of $SU_q(2)$ to build spectral triples for $SU_q(2)/U(1)$ and $U_q(2)$. The algebraic approach leads to a $0$-summable triple in the first case whereas the geometric approach leads to a regular $4$-dimensional theory in the second case. \\

\textbf{3.2. Geometry of $SU_q(2)$.} We use \cite{DV} to construct a Hilbert space $\h$ and a faithful representation $\rho$ of $SU_q(2)$ on $\h$. Let us fix the unitary representations $(V_{q,l}, \pi_{q,l})$ of $U_q(\mathfrak{su_2})$ as in 2.4. In the notation of \cite{DV} the vector space $\mathbb{C}[SU_q(2)]$ is considered in the form 
\begin{eqnarray*}
\mathbb{C}[SU_q(2)] = \bigoplus_{l \in \frac{1}{2}\mathbb{N}_0} V_{q,l} \otimes V^*_{q,l}.
\end{eqnarray*}
The basis is chosen by
\begin{eqnarray*} 
t^l_{m,n} =  | l, m \rangle \otimes  \langle l, n | : \5 l \in \frac{1}{2}\mathbb{N}_0,\5 -l \leq m,n \leq l.
\end{eqnarray*}
and the coproduct is $t^l_{m,n} = \sum_k t^l_{m,k} \otimes t^l_{k,n}$. We notice a conceptual difference compared to our conventions: the second component is treated as a dual meaning that the natural pairing is defined by 
\begin{eqnarray*}
t^l_{m,n}(x) = \langle l,n | \pi_{q,l}(x) | l,m \rangle, 
\end{eqnarray*}
for all $x \in U_q(\mathfrak{su_2})$ and $t^l_{m,n} \in \mathbb{C}[SU_q(2)]$. However, it is straightforward to apply the general theory. We just define a natural left action by
\begin{eqnarray*}
\partial(x)t^l_{m,n} = \pi_{q,l}(x)|l,m \rangle \otimes \langle l,n |,
\end{eqnarray*}
or equivalently, $\partial(x)t = t'(x)t''$ for all $t \in \mathbb{C}[SU_q(2)]$. 

Denote by $C_q$ the unitary Clebsch-Gordan matrices for the representations $\pi_{q,l}$ so that the multiplication is derived from the formulas
\begin{eqnarray*}
t^l_{m,n} t^{l'}_{m',n'} = \sum_{p = |l-l'|}^{l+l'} C_q\begin{pmatrix} 
  l & l' & p \\
  m & m' & m+m' \end{pmatrix}  C_q\begin{pmatrix} 
  l & l' & p \\
  n & n' & n+n' \end{pmatrix} t^{p}_{m+m',n+n'}.
\end{eqnarray*}
The vector $t^0_{00} = 1$ is the unit vector. The inner product of $\mathbb{C}[SU_q(2)]$ given by the Haar state is fixed by
\begin{eqnarray*}
\langle t^{l}_{m,n}, t^{l'}_{m',n'} \rangle:= h((t^l_{m,n})^* t^{l'}_{m',n'})  = \frac{q^{-2m}}{[2l+1]} \delta_{l l'} \delta_{m m'} \delta_{n n'}.
\end{eqnarray*}
Applying the Clebsch Gordan coefficients with the inner product we find the involution
\begin{eqnarray*}
(t^l_{m,n})^* = (-1)^{2l+m+n}q^{n-m}t^l_{-m,-n}.
\end{eqnarray*}
Let us fix $a = t^{\frac{1}{2}}_{\frac{1}{2},\frac{1}{2}}$ and $b = t^{\frac{1}{2}}_{\frac{1}{2},-\frac{1}{2}}$. It follows that $a^* = t^{\frac{1}{2}}_{-\frac{1}{2}, -\frac{1}{2}}$ and $ - q b^* = t^{\frac{1}{2}}_{-\frac{1}{2}, \frac{1}{2}}$ and 
\begin{eqnarray}\label{algsu2}
ba = qab, \5 b^*a = qab^*,\5 b b^* = b^* b,\5 a^* a + q^2 b^* b = 1, \5 a  a^* + b b^* = 1
\end{eqnarray}
determine the algebraic structure. Thus, $a$ and $b$ generate $\mathbb{C}[SU_q(2)]$ as a $*$-algebra.

The orthonormal basis of the prehilbert space $\mathbb{C}[SU_q(2)]$ is
\begin{eqnarray*}
|lmn \rangle = q^m [2l+1]^{\frac{1}{2}} t^l_{m,n}. 
\end{eqnarray*}
Denote by $\h$ the Hilbert space completion of $\mathbb{C}[SU_q(2)] \otimes \Sigma$, where $\Sigma = V_{q,\frac{1}{2}}$. The representation of $U_q(\g)$ on $\h$ is defined by $x \mapsto (\partial \otimes \pi_{q,\frac{1}{2}}) \triangle_q(x)$. The prehilbert space decomposes into irreducible components under this action as
\begin{eqnarray*}
(\bigoplus_{l \in \frac{1}{2}\mathbb{N}_0} V_{q,l} \otimes V^*_{q,l}) \otimes \Sigma \simeq V_{q,\frac{1}{2}} \oplus \bigoplus_{j \in \frac{1}{2} \mathbb{N}} (V_{q,j+\frac{1}{2}} \otimes V^*_{q,j}) \oplus (V_{q,j-\frac{1}{2}} \otimes V^*_{q,j}) := W^{\uparrow}_0 \oplus \bigoplus_{j \in \frac{1}{2}\mathbb{N}} W^{\uparrow}_j \oplus W^{\downarrow}_j.
\end{eqnarray*}
The components $W^{\uparrow}_j$ and $W^{\downarrow}_j$ have multiplicities $(2j+2)(2j+1)$ and $2j(2j+1)$. The orthonormal basis of $\h$ is chosen by
\begin{eqnarray} \label{dec}
 |j \mu n \uparrow \rangle , |j' \mu' n \downarrow \rangle : j \in \frac{1}{2} \mathbb{N}_0,\ j' \in  \frac{1}{2} \mathbb{N},\ |\mu| \leq j+1,\ |\mu'| \leq j-1,\ |n| \leq j 
\end{eqnarray}
where $|j \mu n \uparrow \rangle \in W^{\uparrow}_j$ and $|j \mu n \downarrow \rangle \in W^{\downarrow}_j$. This spectral decomposition was also used in \cite{DV}, where the Dirac operator on $\h$ was fixed from the condition that $[D,\rho(x)]$ is a bounded operator for all $x \in \mathbb{C}[SU_q(2)]$.  

The Dirac operators $\mathfrak{D}_q$ and $D_q$ defined in 2.4. act on $\h$. It is straightforward to compute their spectrum on an arbitrary irreducible component of $\h$
\begin{eqnarray*}
& &\mathfrak{D}_q|j \mu n \uparrow \rangle = [2j]|j \mu n \uparrow \rangle \\
& &\mathfrak{D}_q|j \mu n \downarrow \rangle = [-(2j+2)]|j \mu n \downarrow \rangle.
\end{eqnarray*}
and therefore 
\begin{eqnarray*}
& &D_q|j \mu n \uparrow \rangle = (2j + \frac{3}{2})|j \mu n \uparrow \rangle \\
& &D_q|j \mu n \downarrow \rangle = (-(2j+2) + \frac{3}{2})|j \mu n \downarrow \rangle.
\end{eqnarray*}
The operator $D_q$ is exactly the same Dirac operator which was defined in \cite{DV}. 

In \cite{DV} a faithful $*$-representation was derived from equivariance conditions with the $U_q(\mathfrak{su_2})$ action but it was also noted that the representation coincides with the one coming directly from the GNS construction. Thus it fits into the general theory \cite{NT}. The representaiton has the following form 
\begin{eqnarray}\label{rep2}
& & \rho(a) := \rho(a_+) + \rho(a_-),\5 \rho(a^*) := \rho(a^*_+) + \rho(a^*_-)\\
& & \rho(a_+)|j \mu n \rangle \rangle = \alpha^+_{j \mu n} |j^+ \mu^+ n^+ \rangle \rangle, \5 \rho(a_-)|j \mu n \rangle \rangle = \alpha^-_{j \mu n} |j^- \mu^+ n^+ \rangle \rangle \nonumber \\
& & \rho(b_+)|j \mu n \rangle \rangle = \beta^+_{j \mu n} |j^+ \mu^+ n^- \rangle \rangle, \5 \rho(b_-)|j \mu n \rangle \rangle = \beta^-_{j \mu n} |j^- \mu^+ n^- \rangle \rangle \nonumber \\
& & \rho(b) := \rho(b_+) + \rho(b_-), \5  \rho(b^*) := \rho(b^*_+) + \rho(b^*_-) \nonumber \\
& & \rho(a^*_+)|j \mu n \rangle \rangle = \tilde{\alpha}^+_{j \mu n} |j^+ \mu^- n^- \rangle \rangle, \5 \rho(a^*_-)|j \mu n \rangle \rangle = \tilde{\alpha}^-_{j \mu n} |j^- \mu^- n^- \rangle \rangle \nonumber \\
& & \rho(b^*_+)|j \mu n \rangle \rangle = \tilde{\beta}^+_{j \mu n} |j^+ \mu^- n^+ \rangle \rangle, \5 \rho(b^*_-)|j \mu n \rangle \rangle = \tilde{\beta}^-_{j \mu n} |j^- \mu^- n^+ \rangle \rangle \nonumber \nonumber, \\
& & |j \mu n \rangle \rangle := \begin{pmatrix} 
  |j \mu n \uparrow \rangle \\
  |j \mu n \downarrow \rangle  \end{pmatrix}, \5 j^{\pm} := j \pm \frac{1}{2} \nonumber
\end{eqnarray}
where the matrices $\alpha^{\pm}_{j\mu n}$, $\beta^{\pm}_{j\mu n}$, $\tilde{\alpha}^{\pm}_{j\mu n}$ and $\tilde{\beta}^{\pm}_{j\mu n}$ are defined in \cite{DV} (Proposition 4.4.). We found that the methods of \cite{NT} applied to the isomorphism $\phi$ defined earlier leads to the model \cite{DV}. On the other hand it is known that a different choice of the isomorpshism $\phi$ gives a spectral triple unitarily equivalent to this one.

Let us now turn the attention into Fredholm modules. Since the operators $\mathfrak{D}_q$ and $\tilde{D}_q = D_q - (3/2)\textbf{1}$ have nontrivial kernels we define approximated sign operators by
\begin{eqnarray*}
\mathfrak{F}_q = \frac{\mathfrak{D}_q}{(1 + \mathfrak{D}_q^2)^{\frac{1}{2}}}, \5 F_q = \frac{\tilde{D}_q}{(1 + \tilde{D}_q^2)^{\frac{1}{2}}}.
\end{eqnarray*}
Recall that a Fredholm module $(A, F, \h)$ is called $n$-summable if $n$ is the smallest integer so that the compact operators $F^2-1$ and $[F,\rho(t)]$ are in $L_{n+}(\h)$ for all $t \in A$. \\

\noindent \textbf{Proposition.} The triples $(\mathbb{C}[SU_q(2)], \mathfrak{F}_q, \h)$ and $(\mathbb{C}[SU_q(2)], F_q, \h)$ define $1$- and $3$-summable Fredholm modules and are homotopy equivalent to each other. \\

\noindent Proof. The triple $(\mathbb{C}[SU_q(2)], F_q, \h)$ is a Fredholm module with summability at most $3$ because it is determined by a  $3$-dimensional spectral triple. On the other hand the smallest $n$ for which $F^2-1 \in L_{n+}(\h)$ is $3$. It is shown in \cite{H} that $(\mathbb{C}[SU_q(2)], \mathfrak{F}_q, \h)$ is a $1$-summable Fredholm module. Following the same lines one checks that the family of Fredholm operators $[0,1] \rightarrow B(\h)$ defined by
\begin{eqnarray*}
& &0 \mapsto F_q, \\
& &t \mapsto \frac{[\widetilde{D}_q]_t}{(1+[\widetilde{D}_q]_t^2)^{\frac{1}{2}}}, \5 [\widetilde{D}_q]_t = \frac{q^{t\widetilde{D}_q} - q^{-t\widetilde{D}_q}}{q^t - q^{-t}},  \5 \hbox{for $t \in (0,1]$ }
\end{eqnarray*}
together with $\rho, \mathbb{C}[SU_q(2)]$ and $\h$ defines a familily of Fredholm modules and connects the operators $F_q$ and $\mathfrak{F}_q$. $\5 \square$ \\

By the Proposition and discussion above we have found the explicit relationship between the models \cite{BK, DV, H, NT}. Especially, from the point of view of index theory they all describe the same element in the $K$-homology. \\

\textbf{3.3. Geometry of $S^2_q$.} The standard Podles sphere $\mathbb{C}[S^2_q]$ is the fixed point algebra $\mathbb{C}[SU_q(2)]^{U(1)}$ under the left action of the group $U(1)$. The action on the generators is given by 
\begin{eqnarray*}
a \mapsto e^{i \phi}a,\5 a^* \mapsto e^{-i \phi}a,\5 b \mapsto e^{i \phi}b,\5 b^* \mapsto e^{-i \phi} b^*.  
\end{eqnarray*}
Equivalently we can consider $S^2_q$ as $U_q(\mathfrak{h})$-invariant subalgebra 
\begin{eqnarray*}
\mathbb{C}[S^2_q] = \{t \in \mathbb{C}[SU_q(2)]: \partial(k)t = t \}.
\end{eqnarray*}
Therefore we can choose the generators of $\mathbb{C}[S^2_q]$ by 
\begin{eqnarray*}
A = ab^*, \5 A^* = ba^*, \5 B = B^* = b b^*,
\end{eqnarray*}
which satisfy the algebraic relations  
\begin{eqnarray*}
AB = q^{-2}BA, \5 A^*B = q^2 B A^*,\5 A A^* = q^{-2}B(1 - B),\5 A^*A =B(1-q^2B).
\end{eqnarray*}

The Hilbert space of the theory is the completion of the invariant subspace of $\mathbb{C}[SU_q(2)] \otimes \Sigma$ under the left $U_q(\mathfrak{h})$-action
\begin{eqnarray*}
 (\mathbb{C}[SU_q(2)] \otimes \Sigma)^{U_q(\mathfrak{h})} = \{\Psi: \partial(k) \otimes \mathrm{id} \otimes \pi_{q,\frac{1}{2}}(k) \Psi = \Psi \}.
\end{eqnarray*}
The completion is done with the state which is the restriction of the Haar state on the invariant subspace. Let us denote by $\h^{\mathfrak{h}}$ the Hilbert space. The basis is chosen by
\begin{eqnarray*}
 |ln + \rangle := |l,-\frac{1}{2}, n \rangle \otimes |\frac{1}{2}, \frac{1}{2} \rangle,\5  |ln - \rangle := |l, \frac{1}{2}, n  \rangle \otimes |\frac{1}{2}, - \frac{1}{2} \rangle
\end{eqnarray*}
for all $l \in \mathbb{N}_0$ and $-l \leq n \leq l$. The representation $\rho$ of $\mathbb{C}[SU_q(2)]$ restricts to a faithful $*$-representation of the subalgebra $\mathbb{C}[S^2_q]$ on $\h^{\mathfrak{h}}$. 

We apply the algebraic Dirac operator model here and use the conventions of Section 2.4. Define $\Omega' = p\Omega$, where $p$ is a projection onto the subspace $V_{q,\rho}/V_{q,0}$ where $V_{q,0}$ is the one dimensional subspace of weight zero vectors in $V_{q,\rho}$. Then we define 
\begin{eqnarray*}
\mathfrak{D}_q = (\theta \otimes \sigma) \Omega' = [2] \partial \begin{pmatrix} 
 0 & q^{-\frac{1}{2}} k^{-1}f\\ 
 q^{\frac{1}{2}} k^{-1} e & 0 \end{pmatrix} = [2] \partial \begin{pmatrix} 
 0 & f\\ 
 e & 0
\end{pmatrix}
\end{eqnarray*}
The latter equality can be checked by applying $\mathfrak{D}_q$ on the basis vectors of $\h^{\mathfrak{h}}$. The operator we have found coincides with the one defined in \cite{S}. The triple ($\mathbb{C}[S^2_q], \mathfrak{D}_q, \h^{\mathfrak{h}}$) is a $0$-summabls spectral triple. The chirality operator is defined by $\gamma = \mathrm{Diag}(1,-1)$. \\

\textbf{3.4. Geometry of $U_q(2)$.} Recall that the irreducible representations are parametrized by $(l,c) \in P_+$. The vector space $\mathbb{C}[U_q(2)]$ is spanned by 
\begin{eqnarray*}
t^{l,c}_{m,n} = |l,c, n \rangle \otimes \langle l,c,m | \in V_{q,(l,c)} \otimes V^*_{q,(l,c)}\end{eqnarray*}
where $(l,c) \in P_+$. The product is determined by the formulas 
\begin{eqnarray*}
t^{l,c}_{m,n}t^{l',c'}_{m',n'} =  \sum_{p = |l-l'|}^{l+l'} C_q\begin{pmatrix} 
  l & l' & p \\
  m & m' & m + m' \end{pmatrix}  C_q \begin{pmatrix} 
  l & l' & p \\
  n & n' & n+n' \end{pmatrix} t^{p,c+c'}_{m+m', n+n'},
\end{eqnarray*}
where the $C_q$ are the Clebsch-Gordan coefficients of $U_q(\mathfrak{su_2})$. The unit is $1 = t^{0,0}_{0,0}$. As a $*$-algebra $\mathbb{C}[U_q(2)]$ is generated by
\begin{eqnarray*}
a = t^{\frac{1}{2},\frac{1}{2} }_{\frac{1}{2}, {\frac{1}{2}}}, \5 b = t^{\frac{1}{2},\frac{1}{2}}_{\frac{1}{2}, -{\frac{1}{2}}}, \5 C = t^{0,1}_{0,0}.
\end{eqnarray*}
$a$ and $b$ satisfy \eqref{algsu2} and $C$ has the properties 
\begin{eqnarray*}
C t^{l,c}_{m,n}  = t^{l,c+1}_{m,n},\5 C^* = t^{0,-1}_{0,0},\5 Ct = tC,
\end{eqnarray*} 
for all $t \in \mathbb{C}[U_q(2)]$.

Denote by $h$ the Haar state of $\mathbb{C}[SU_q(2)]$. We extend this to the Haar state of $\mathbb{C}[U_q(2)]$ by
\begin{eqnarray*}
\hat{h}(t^{l,c}_{m,n}) = \delta_{c,0} h(t^{l}_{m, n}).
\end{eqnarray*}
It is left invarint $(\hat{h} \otimes \mathrm{id})(\triangle_q(t)) = \hat{h}(t) 1$ for all $t \in \mathbb{C}_q[U(2)]$. The involution is given by $(t^{l,c}_{m,n})^* = (-1)^{2l+m+n}q^{n-m}t^{l,-c}_{-m,-n}$. Furthermore, 
\begin{eqnarray*}
||t^{l,c}_{m,n}||^2 = \hat{h}((t^{l,c}_{m,n})^* t^{l,c}_{m,n}) = h((t^{l}_{m,n})^* t^{l}_{m,n}), 
\end{eqnarray*}
which vanishifs if and only if $t^{l}_{m,n} = 0$. Thus, the Haar state $\hat{h}$ is faithful. Let $L^2_q(U(n))$ be the Hilbert space completion. The orthonormal basis is defined by
\begin{eqnarray*}
|lmnc \rangle = q^m [2l+1]^{\frac{1}{2}} t^{l,c}_{m,n}.
\end{eqnarray*}

Recall that the spinor module is $\hat{\Sigma} \simeq V^+_{q,(\frac{1}{2},0)} \oplus V^-_{q,(\frac{1}{2},0)}$. We define subspaces $W^{\uparrow}_{j,c,\pm}$ and $W^{\downarrow}_{j,c,\pm}$ so that the decomposition onto irreducible components under the left action is  
\begin{eqnarray*}
&&(V_{q,(j,c)} \otimes V^*_{q,(j,c)}) \otimes V^{\pm}_{q,(\frac{1}{2},0)} \simeq  (V^{\pm}_{q,(j+\frac{1}{2},c)} \otimes V^*_{q,(j,c)}) \oplus ( V^{\pm}_{q,(j-\frac{1}{2},c)} \otimes V^*_{q,(j,c)}) \\
&&= W^{\uparrow}_{j,c,\pm} \oplus  W^{\downarrow}_{j,c,\pm}
\end{eqnarray*}
and then the Hilbert space decomposes into irreducible compenents by
\begin{eqnarray*}
\widehat{\h} = L^2(U_q(2)) \otimes \hat{\Sigma} = \bigoplus_{c}\Big( W^{\uparrow}_{0,c,+} \oplus  W^{\uparrow}_{0,c,-} \oplus \bigoplus_{j \in \frac{1}{2}\mathbb{N}}^{\infty} W^{\uparrow}_{j,c,+} \oplus W^{\downarrow}_{j,c,+} \oplus W^{\uparrow}_{j,c,-} \oplus W^{\downarrow}_{j,c,-} \Big).
\end{eqnarray*}
Now the sum over $c$ is defined so that $c$ runs over half integers for each integer $j$ and over integers for each half integer $j$. We can fix an orthonormal basis 
\begin{eqnarray*}
|j \mu n \uparrow c \pm  \rangle, |j' \mu' n  \downarrow c \pm \rangle :\5 j,j', \mu, \mu',n, c 
\end{eqnarray*}
so that $j,j', \mu,\mu'$ and $n$ are restricted as in \eqref{dec}. Let us adopt a column vector notation
\begin{eqnarray*}
|j \mu n c \rangle \rangle = \begin{pmatrix} 
  |j \mu n \uparrow c + \rangle \\
  |j \mu n \downarrow c + \rangle \\
  |j \mu n \uparrow c - \rangle \\
  |j \mu n \downarrow c - \rangle \\  
  \end{pmatrix} = \begin{pmatrix} 
  |j \mu n c + \rangle \rangle \\
  |j \mu n c - \rangle \rangle
  \end{pmatrix}
\end{eqnarray*}
The algebra of functions has a diagonal action on $\h$
\begin{eqnarray*}
\hat{\rho}(t)|j \mu n c\rangle \rangle = \begin{pmatrix} 
  \rho(t) & 0 \\
   0 & \rho(t) \end{pmatrix} |j \mu n, c + \frac{1}{2} \rangle \rangle,\5 \hat{\rho}(t^*)|j \mu n c\rangle \rangle = \begin{pmatrix} 
  \rho(t^*) & 0 \\
   0 & \rho(t^*) \end{pmatrix} |j \mu n, c - \frac{1}{2} \rangle \rangle
\end{eqnarray*}
where for $t =a,b$ and the representation $\rho(t)$ is independent of the parameters $c, \pm$ and given as in \eqref{rep2} whereas the generator $C$ acts by
\begin{eqnarray*}
\hat{\rho}(C)|j \mu n c\rangle \rangle = |j \mu n, c+1 \rangle \rangle, \5 \hat{\rho}(C^*)|j \mu n c\rangle \rangle = |j \mu n, c-1 \rangle \rangle.
\end{eqnarray*}
Recall that the geometric Dirac operator $\hat{D}_q$ has the following form
\begin{eqnarray*}
\hat{D}_q |j \mu n c \rangle \rangle =  \begin{pmatrix} 
  0 & iD_q + \partial(x_0) \\
  -iD_q + \partial(x_0) & 0  \end{pmatrix} |j \mu n c \rangle \rangle
\end{eqnarray*}
where $D_q$ is the Dirac operator on $SU_q(2)$. The absolute value operator and the chirality are defined by
\begin{eqnarray*}
& & |\hat{D}_q| = \begin{pmatrix} 
  (D_q^2 + \partial(x_0)^2)^{\frac{1}{2}} & 0 \\
  0 & (D_q^2 + \partial(x_0)^2)^{\frac{1}{2}} \end{pmatrix}, \5  \gamma  = \begin{pmatrix} 
  \textbf{1} & 0 \\
  0 & -\textbf{1} \end{pmatrix} 
\end{eqnarray*}
The chirality satisfies 
$\gamma \hat{D}_q + \hat{D}_q \gamma = 0,\ \gamma^2 = 1,\ \gamma = \gamma^*$ and $[\gamma, \hat{\rho}(t)] = 0$ as it should in a $4$-dimensional model. 

We define following projection operators
\begin{eqnarray*}
\hat{P}^{\uparrow} = \begin{pmatrix} 
  P^{\uparrow} & 0 \\
   0 & P^{\uparrow} \end{pmatrix},\5 \hat{P}^{\downarrow} = \begin{pmatrix} 
  P^{\downarrow} & 0 \\
   0 & P^{\downarrow} \end{pmatrix}
\end{eqnarray*}
where $P^{\uparrow}$ and $P^{\downarrow}$ are projection operators onto positive and negative eigenspaces for $D_q$. \\

\noindent \textbf{Proposition.} The triple $(\mathbb{C}[U_q(2)],\hat{D}_q, \gamma, \h)$ is an even and regular $4$-dimensional spectral triple. \\

\noindent Proof. The dimensionality is true by construction. We need to check that $[\hat{D}_q, \hat{\rho}(t)]$ is bounded for each $t \in \mathbb{C}[U_q(2)]$. It is sufficient to prove this for the generators. We have
\begin{eqnarray*}
[\hat{D}_q,\hat{\rho}(t)] =  \begin{pmatrix} 
  0 & i[D_q, \rho(t)] \\
  -i[D_q,\rho(t)] & 0  \end{pmatrix} + \begin{pmatrix} 
  0 & [\partial(x_0), \rho(t)] \\
  [\partial(x_0),\rho(t)] & 0  \end{pmatrix}
\end{eqnarray*}
and the first term is bounded by the $SU_q(2)$ theory and the second term is bounded because $[\partial(x_0), \rho(t)] = \frac{1}{2} \rho(t)$ if $t = a,b$. Furthremore, 
\begin{eqnarray*}
[\hat{D}_q, \hat{\rho}(C)] = \begin{pmatrix} 
  0 & 1 \\
   1 & 0  \end{pmatrix} \hat{\rho}(C)
\end{eqnarray*}
is bounded. 

To prove the regularity we need to show that the operators $\delta^k(\hat{\rho}(t))$ and $\delta^k([\hat{D}_q, \hat{\rho}(t)])$ are bounded for each $k \in \mathbb{N}$. Since $\delta$ is a derivation it is enough to check this for the generators. It is known that the off diagonal operators (wrt. the polarization into positive and negative eigenspaces) $P^{\uparrow} \rho(t) P^{\downarrow}$ and $P^{\downarrow} \rho(t) P^{\uparrow}$ are given by rapidly decaying sequances for any $t = a, b$. Consequently $\delta^p(P^{\uparrow} \rho(t) P^{\downarrow})$ and $\delta^p(P^{\downarrow} \rho(t) P^{\uparrow})$
are bounded and even trace class because of the polynomial growth of the eigenvalues of $D_q$. Let $t_{\pm} = a_{\pm}$ or $b_{\pm}$, recall \eqref{rep2}. Then 
\begin{eqnarray*}
& &\delta^p(\hat{\rho}(t_{+}))|j \mu n c \rangle \rangle \approx \Big[(n_{j+\frac{1}{2},c + \frac{1}{2}} - n_{j,c})^p P^{\uparrow} \hat{\rho}(t_+) P^{\uparrow} + (n_{j,c + \frac{1}{2}} - n_{j-\frac{1}{2},c})^p  P^{\downarrow} \hat{\rho}(t_+) P^{\downarrow} \Big] |j \mu n c \rangle \rangle, \\
& &\delta^p(\hat{\rho}(t_{-})) |j \mu n c \rangle \rangle \approx \Big[( n_{j-\frac{1}{2},c + \frac{1}{2}} - n_{j,c})^p P^{\uparrow} \hat{\rho}(t_-) P^{\uparrow} + (n_{j-1,c + \frac{1}{2}} - n_{j-\frac{1}{2},c})^p P^{\downarrow} \hat{\rho}(t_-) P^{\downarrow} \Big]|j \mu n c \rangle \rangle, \nonumber\\
& &\delta^p( \hat{\rho}(C))|j \mu n c \rangle \rangle \approx \Big[(n_{j,c+1} - n_{j,c})^p P^{\uparrow} \hat{\rho}(C) P^{\uparrow} + (n_{j-\frac{1}{2},c+1} - n_{j-\frac{1}{2},c})^p P^{\downarrow} \hat{\rho}(C) P^{\downarrow} \Big] |j \mu n c \rangle \rangle \nonumber
\end{eqnarray*}
where the symbol $\approx$ means that the equality holds modulo trace class contributions. The constants are defined by
\begin{eqnarray*}
n_{j,c} = \sqrt{(2j+3/2)^2 + c^2}.
\end{eqnarray*}
The functions such as $n_{j+\frac{1}{2},c + \frac{1}{2}} - n_{j,c}$ above are bounded for all values of $j,c$. The operators $P^{\uparrow} \hat{\rho}(t) P^{\uparrow}$ and $P^{\downarrow} \hat{\rho}(t) P^{\downarrow}$ in the formulas above are bounded by the $SU_q(2)$ theory and therefore $\delta^p(\hat{\rho}(t))$ is bounded for any $p$ and $t = a,b, C$. Finally, applying $\delta([\hat{D}_q,\hat{\rho}(t)] = [\hat{D}_q,\delta(\hat{\rho}(t))]$ and consequently  $\delta^p([\hat{D}_q,\hat{\rho}(t)] = [\hat{D}_q,\delta^p(\hat{\rho}(t))]$ with the above operators one can immediately check that $\delta^p([\hat{D}_q,\hat{\rho}(t)]$ is bounded for any $t = a,b, C$. Therefore the spectral triple is regular. $\5 \square$ \\

\textbf{Acknowledgements.} This project was completed during the first half of 2011 at Cardiff University and is funded by the Marie Curie Research Training Network MRTN-CT-2006-031962 in Noncommutative Geometry (EU-NCG). The author wishes to thank Prof. David E. Evans and Dr. Mathew Pugh for support during the project.


\begin{thebibliography}{1}
\bibitem{BK} Bibikov P.N., Kulish P.P.: Dirac operators on quantum SU(2) group and quantum sphere, J. Math. Sci, 100, 2039-2050 (2000)
\bibitem{CP} Chakraborty P.S., Pal A.: Equivariant Spectral Triples on the Quantum SU(2) Group, K-Theory 28, 107-126 (2003)
\bibitem{CP2} Chakraborty P.S., Pal A.: On equivariant Dirac operators for $SU_q(2)$, Proc. Indian Acad. Sci. 116, 531-541 (2006)
\bibitem{Con} Connes A.: Gravity Coupled with Matter and the Foundation of Non-commutative Geometry, Commun. Math. Phys. 182, 155-176 (1996)
\bibitem{Cu} Curtright T. L., Ghandour G. I., Zachos G. K.: Quantum algebra deforming maps, Clebsch-Gordan coefficients, coproducts, R and U matrices, J. Math. Phys. 32, 676-688 (1991)
\bibitem{DV} Dabrowski L., Landi G., Sitarz A., van Suijlekom W., Varilly J.C.: The Dirac Operator on $SU_q(2)$, Commun. Math. Phys. 259, 729-759 (2005)
\bibitem{S}  Dabrowski L, Sitarz, A.: Dirac operator on the standard Podles quantum sphere, Noncommutative geometry and quantum groups (Warsaw, 2001), Banach Center Publ. 61, 49-58 (Polish Acad. Sci. 2003)
\bibitem{DM} Deligne P, Milne J.S.: Tannakian Categories, in: Hodge Cycles, Motives and Shimura Varieties, Letcture Notes in Mathematics 900, 101-228 (Springer-Verlag 1981)
\bibitem{De} Delius G.W., Gould M.D.: Quantum Lie Algebras, Their Existence, Uniqueness and $q$-antisymmetry, Commun. Math. Phys. 185, 709-722 (1997)
\bibitem{Dr} Drinfeld V. G.: Quantum groups. Proc. I.C.M. Berkeley (1986)
\bibitem{H} Harju A.: Covariant Dirac Operators on Quantum Groups, arXiv:1009.3913
\bibitem{HP} Huang J-S., Pandzic P.: Dirac Operators in Representation Theory (Birkhauser 2006)
\bibitem{Ji} Jimbo M.: A $q$-Difference Analogue of $U(g)$ and the Yang-Baxter Equation, Lett. Math. Phys. 10, 63-69 (1985)
\bibitem{KR} Kirillov A. N., Reshetikhin N.: $q$-Weyl group and a multiplicative formula for universal $R$-matrices, Commun. Math. Phys. 134, 421-431 (1990)
\bibitem{Lu} Lusztig G.: Quantum Deformations os Certain Simple Modules Over Enveloping Algebras, Adv. Math. 70, 237-249 (1988)
\bibitem{KL1} Kazhdan D., Lusztig G.: Tensor Structures Arising from Affine Lie Algebras. III, J. Amer. Math. Soc. 7, 335-381 (1994)
\bibitem{KL2} Kazhdan D., Lusztig G.: Tensor Structures Arising from Affine Lie Algebras. IV, J. Amer. Math. Soc. 7, 383-453 (1994)
\bibitem{KS} Korogodski L. I., Soibelman Y. S.: Algebras of Functions on Quantum Groups: Part I (AMS 1998)
\bibitem{N} Nagy G.: Deformation quantization and K-theory, in: Perspectives on quantization (South Hadley, MA, 1996),
Contemp. Math. 214, 111-134 (AMS 1998)
\bibitem{NT0} Neshveyev S., Tuset L.: Notes on the Kazhdan-Lusztig theorem on equivalence of the Drinfeld category and the category of $U_q(\g)$-modules, arXiv:0711.4302v1
\bibitem{NT} Neshveyev S., Tuset L.: The Dirac operator on compact quantum groups, J. Reine Angew. Math. 641, 1-20 (2010)
\bibitem{NT2} Neshveyev S., Tuset L.: Symmetric Invariant Cocycles on the Duals of $q$-deformations, Adv. Math. 227, 146-169 (2011)
\bibitem{NT3} Neshveyev S., Tuset L.: K-homology class of the Dirac operator on a compact quantum group, arXiv:1102.0248
\bibitem{Ro} Rosso M.: Finite Dimensional Representations of the Quantum Analog of the Enveloping Algebra of a Complex Simple Lie Algebra, Commun. Math. Phys. 117, 581-593 (1988)
\end{thebibliography}
\end{document}